\documentclass[a4paper,12pt]{article}
\usepackage{amsfonts,amssymb,latexsym,amsmath}
\usepackage{multicol}
\usepackage[cp1251]{inputenc}
\usepackage[russian]{babel}
\begin{document}
\renewcommand{\refname}{References}

\date{}
\title{The invariant field of the adjoint action of the
unitriangular group in the nilradical of a~parabolic subalgebra}
\author{V. V. Sevostyanova\thanks{This research was partially
supported by the RFBR (project 08-01-00151-а)}}

\maketitle

\begin{center}
\parbox[b]{330pt}{\small\textsc{Abstract.} In the present paper
the field of invariants of the adjoint action of the unitriangular
group in the nilradical of any parabolic subalgebra is described.}
\end{center}

\vspace{0.5cm}

Consider the general linear group $\mathrm{GL}(n,K)$ defined over an
algebraically closed field $K$ of characteristic 0. Let $B$ ($N$,
respectively) be its Borel (maximal unipotent, respectively)
subgroup, which consists of triangular matrices with nonzero (unit,
respectively) elements on the diagonal. We fix a parabolic subgroup
$P$ that contains $B$. Denote by $\mathfrak{p}$, $\mathfrak{b}$ and
$\mathfrak{n}$ the Lie subalgebras in $\mathfrak{gl}(n,K)$ that
correspond to $P$, $B$ and $N$ respectively. We represent
$\mathfrak{p}=\mathfrak{r}\oplus\mathfrak{m}$ as the direct sum of
the nilradical $\mathfrak{m}$ and a block diagonal subalgebra
$\mathfrak{r}$ with sizes of blocks $(n_1,\ldots, n_s)$. The
subalgebra $\mathfrak{m}$ is invariant relative to the adjoint
action of the group $P$, therefore, $\mathfrak{m}$ is invariant
relative to the action of the subgroups $B$ and $N$. We extend this
action to the representation in the algebra $K[\mathfrak{m}]$ and in
the field $K(\mathfrak{m})$. The subalgebra $\mathfrak{m}$ contains
a Zariski-open $P$-orbit, which is called the \emph{Richardson
orbit} (see~\cite{R}). Consequently, the algebra of invariants
$K[\mathfrak{m}]^P$ coincides with $K$. The question concerning the
structure of the algebra of invariantss $K[\mathfrak{m}]^N$ and
$K[\mathfrak{m}]^B$ remains open and seems to be a considerable
challenge. In the special case $P=B$, the answer is following: the
algebra of invariants $K[\mathfrak{m}]^N$ is the polynomial algebra
$K[x_{12},x_{23},\ldots,x_{n-1,n}]$. The present paper is a
continuation of a paper~\cite{PS}, in which a conjecture on the
structure of the field of invariants $K(\mathfrak{m})^N$ was stated.
This conjecture is proved in the present paper. Theorem 1.8 provides
a complete description of the field of invariants
$K(\mathfrak{m})^N$ for any parabolic subalgebra. Another result of
the paper is Theorem 1.7, in which canonical representation of
$N$-orbits in general position are indicated.

\section*{1. The main definitions and results}

\hspace{1.2em} Every positive root $\gamma$ in $\mathfrak{gl}(n,K)$
has the form~(see \cite{GG}) $\gamma=\varepsilon_i-\varepsilon_j$,
$1\leqslant i<j\leqslant n$. We identity a root $\gamma$ with the
pair $(i,j)$ and the set of the positive roots $\Delta^{\!+}$ with
the set of pairs $(i,j)$, $i<j$. The system of positive roots
$\Delta^{\!+}_\mathfrak{r}$ of the reductive subalgebra
$\mathfrak{r}$ is a subsystem in $\Delta^{\!+}$.

Let $\{E_{i,j}:~i<j\}$ be the standard basis in $\mathfrak{n}$. By
$E_\gamma$ denote the basis element $E_{i,j}$, where $\gamma=(i,j)$.

We define a relation in $\Delta^{\!+}$ such that
$\gamma'\succ\gamma$ whenever
$\gamma'-\gamma\in\Delta^{\!+}_\mathfrak{r}$. If
$\gamma\prec\gamma'$ or $\gamma\succ\gamma'$, then the roots
$\gamma$ and $\gamma'$ \emph{comparable}. Denote by $M$ the set of
$\gamma\in\Delta^{\!+}$ such that $E_\gamma\in\mathfrak{m}$. We
identify the algebra $K[\mathfrak{m}]$ with the polynomial algebra
in the variables $x_{i,j}$, ~$(i,j)\in M$.

\medskip
\textbf{Definition 1.1.} A subset $S$ in $M$ is called a \emph{base}
if the elements in $S$ are not pairwise comparable and for any
$\gamma\in M\setminus S$ there exists $\xi\in S$ such that
$\gamma\succ\xi$.

\medskip
\textbf{Definition 1.2.} Let $A$ be a subset in $S$. We say that
$\gamma$ is a \emph{minimal element} in $A$ if there is no $\xi\in
A$ such that $\gamma\succ\xi$.

\medskip
Note that $M$ has a unique base $S$, which can be constructed in the
following way. We form the set $S_1$ of minimal elements in $M$. By
definition, $S_1\subset S$. Then we form a set $M_1$, which is
obtained from $M$ by deleting $S_1$ and all elements
$$\{\gamma\in M:\exists\ \xi\in S_1,\ \gamma\succ\xi\}.$$
The set of minimal elements $S_2$ in $M_1$ is also contained in $S$,
and so on. Continuing the process, we get the base $S$.

\medskip
\textbf{Definition 1.3.} An ordered set of positive roots
$\{\gamma_1,\ldots,\gamma_s\}$ is called a \emph{chain} if
$\gamma_1=(a_1,a_2)$, $\gamma_2=(a_2,a_3)$, $\gamma_3=(a_3,a_4)$,
and so on. The number $s$ is called the \emph{length of a chain}.

\medskip
\textbf{Definition 1.4.} We say that two roots $\xi,\xi'\in S$ form
an \emph{admissible pair} $q=(\xi,\xi')$ if there exists
$\alpha_q\in\Delta^{\!+}_\mathfrak{r}$ such that the ordered set of
roots $\{\xi,\alpha_q,\xi'\}$ is a chain. Note that the root
$\alpha_q$ is uniquely determined by $q$.

\medskip
We form the set $Q:=Q(\mathfrak{p})$ that consists of admissible
pairs of roots in $S$. For every admissible pair $q=(\xi,\xi')$ we
construct a positive root $\phi_q=\alpha_q+\xi'$. Consider the
subset $\Phi=\{\phi_q:~ q\in Q\}$.

Using the given parabolic subgroup, we construct a diagram, which is
a square matrix in which the roots from $S$ are marked by the symbol
$\otimes$ and the roots from $\Phi$ are labeled by the symbol
$\times$. The other entries in the diagram are empty.

\textbf{Example 1.} Below a diagram for a parabolic subalgebra with
sizes of its diagonal blocks $(1,3,2,1,3,2,2)$ is given.
\begin{center}
{\begin{tabular}{|p{0.1cm}|p{0.1cm}|p{0.1cm}|p{0.1cm}|p{0.1cm}|
p{0.1cm}|p{0.1cm}|p{0.1cm}|p{0.1cm}|p{0.1cm}|p{0.1cm}|p{0.1cm}|
p{0.1cm}|p{0.1cm}|c} \multicolumn{2}{l}{{\small 1\hspace{5pt}
2\!\!}}&\multicolumn{2}{l}{{\small 3\hspace{5pt}
4\!\!}}&\multicolumn{2}{l}{{\small 5\hspace{5pt} 6\!\!}}&
\multicolumn{2}{l}{{\small 7\hspace{5pt}
8\!\!}}&\multicolumn{2}{c}{{\small\makebox[0.2cm][c]{9\ \
10\!}}}&\multicolumn{2}{c}{{\small\makebox[0.2cm][c]{11\hspace{3pt}12}}}&
\multicolumn{2}{c}{{\small\makebox[0.2cm][c]{13\hspace{3pt}14}}}\\
\cline{1-14} 1&$\otimes$&&&&&&&&&&&&&{\small 1}\\
\cline{1-14} &\multicolumn{3}{l|}{1}&$\times$&$\times$&&&&$\otimes$&
&&&&{\small 2}\\
\cline{5-14} &\multicolumn{3}{c|}{1}&&$\otimes$&&&&&&&&&{\small 3}\\
\cline{5-14} &\multicolumn{3}{r|}{1}&$\otimes$&&&&&&&&&&{\small 4}\\
\cline{2-14} \multicolumn{4}{|c|}{}&\multicolumn{2}{l|}{1}&$\times$&
&$\otimes$&&&&&&{\small 5}\\
\cline{7-14}
\multicolumn{4}{|c|}{}&\multicolumn{2}{r|}{1}&$\otimes$&
&&&&&&&{\small 6}\\
\cline{5-14} \multicolumn{6}{|c|}{}&1&$\otimes$&&&&&&&{\small 7}\\
\cline{7-14} \multicolumn{7}{|c|}{}&\multicolumn{3}{l|}{1}&$\times$&
$\times$&&&{\small 8}\\
\cline{11-14}
\multicolumn{7}{|c|}{}&\multicolumn{3}{c|}{1}&$\times$&
$\otimes$&&&{\small 9}\\
\cline{11-14}
\multicolumn{7}{|c|}{}&\multicolumn{3}{r|}{1}&$\otimes$&
&&&{\small 10}\\
\cline{8-14}
\multicolumn{10}{|c|}{}&\multicolumn{2}{l|}{1}&$\times$&
$\otimes$&{\small 11}\\
\cline{13-14} \multicolumn{10}{|c|}{}&\multicolumn{2}{r|}{1}&
$\otimes$&&{\small 12}\\
\cline{11-14} \multicolumn{12}{|c|}{}&\multicolumn{2}{l|}{1}&
{\small 13}\\
\multicolumn{12}{|c|}{}&\multicolumn{2}{r|}{1}&{\small 14}\\
\cline{1-14} \multicolumn{14}{c}{Diagram 1.}\\
\end{tabular}}
\end{center}

Consider the formal matrix $\mathbb{X}$ in which the variables
$x_{i,j}$ occupy the positions $(i,j)\in M$ and the other entries
are equal to zero. For any root $\gamma=(a,b)\in M$ we denote by
$S_\gamma$ the set of $\xi=(i,j)\in S$ such that $i>a$ и $j<b$. Let
$S_\gamma=\{(i_1,j_1),\ldots,(i_k,j_k)\}$. Denote by $M_\gamma$ a
minor $M_I^J$ of the matrix $\mathbb{X}$ with ordered systems of
rows $I$ and columns $J$, where
$$I=\mathrm{ord}\{a,i_1,\ldots,i_k\},\quad J=\mathrm{ord}\{j_1,\ldots,j_k, b\}.$$

For every admissible pair $q=(\xi,\xi')$, we construct the
polynomial
\begin{equation}
L_q=\sum_{\scriptstyle\alpha_1,\alpha_2\in\Delta^{\!+}_\mathfrak{r}\cup\{0\}
\atop\scriptstyle\alpha_1+\alpha_2=\alpha_q}
M_{\xi+\alpha_1}M_{\alpha_2+\xi'}.\label{L_q}
\end{equation}

\medskip
\textbf{Theorem 1.5}~\cite{PS}\textbf{.} \emph{For an arbitrary
parabolic subalgebra}, \emph{the system of polynomials}
$$\{M_\xi,~\xi\in S,~L_{\varphi},~\varphi\in\Phi,\}$$
\emph{is contained in $K[\mathfrak{m}]^N$ and is algebraically
independent over $K$}.

\medskip
Denote by $\mathcal{Y}$ the subset in $\mathfrak{m}$ that consists
of matrices of the form
$$\sum_{\xi\in S}c_{\xi}E_\xi+
\sum_{\varphi\in\Phi}c'_{\varphi}E_\varphi,$$ where $c_{\xi}\neq0$
and $c'_{\varphi}\neq0$.

\medskip
\textbf{Definition 1.6.} The matrices from $\mathcal{Y}$ are said to
be \emph{canonical}.

\medskip
In the fourth section, we prove the following two theorems.

\medskip
\textbf{Theorem 1.7.} \emph{There exists a nonempty Zariski-open
subset $U\subset\mathfrak{m}$ such that the $N$-orbit of any $x\in
U$ intersects $\mathcal{Y}$ at a unique point}.

\medskip
\textbf{Theorem 1.8.} \emph{The field of invariants
$K(\mathfrak{m})^N$ is the field of rational functions of} $M_\xi$,
$\xi\in S$, \emph{and} $L_{\varphi}$, $\varphi\in\Phi$.

\medskip
The following result is a consequence of Theorem 1.8.

\medskip
\textbf{Theorem 1.9.} \emph{The maximal dimension of an $N$-orbit in
$\mathfrak{m}$ is equal to} $\mathrm{dim}\mathfrak{m}-|S|-|\Phi|$.

\medskip
We introduce on the set of roots $S\cup\Phi$ an order relation for
which
\begin{itemize}
\item[1)] $\xi<\varphi$ for any $\xi\in S$ and $\varphi\in\Phi$;
\item[2)] for other pairs of roots from $S\cup \Phi$, the relation
$<$ means that $(a,b)<(c,d)$ if $c<a$ or $c=a$ if $b<d$.
\end{itemize}

Consider the roots from the set $S\cup\Phi$:
$S=\{\xi_1,\ldots,\xi_p\}$, $\Phi=\{\varphi_1,\ldots,\varphi_q\}$.
Let the mapping $\varpi:\mathfrak{m}\rightarrow K^{p+q}$ be given by
$$x\mapsto\Big(M_{\xi_1}(x),\ldots,M_{\xi_p}(x),
L_{\varphi_1}(x),\ldots,L_{\varphi_q}(x)\Big).$$ Denote
$(c_{\xi_1},\ldots,c_{\xi_p},c_{\varphi_1},\ldots,c_{\varphi_q})\in
\varpi(\mathfrak{m})$. Let $\alpha\in S$. By $\widehat{c}_{\alpha}$
denote the number $\prod c_{\beta}$, where the product proceeds over
all roots $\beta\in S$ such that $\beta<\alpha$ and $\beta$ is the
greatest root in the sense of the order~$<$. Let $U$ be a
Zariski-open set in the condition of Theorem 1.7. We have a
canonical representative of the $N$-orbit of any matrix from~$U$.

\medskip
\textbf{Corollary 1.10.} \emph{Let $x\in U$ and}
$\varpi(x)=(c_{\xi_1},\ldots,c_{\xi_p},c_{\varphi_1},\ldots,c_{\varphi_q})$.
\emph{Then the canonical representative of the $N$-orbit of the
element~$x$ is the following element}
$$\sum_{i=1}^{p}\frac{c_{\xi_i}}{\widehat{c}_{\xi_i}}E_{\xi_i}+
\sum_{j=1}^{q}\frac{c_{\varphi_j}}{c_{\gamma_j}\widehat{c}_{\gamma'_j}}E_{\varphi_j},$$
\emph{where $\varphi_j\in\Phi$ corresponds to the admissible pair}
$(\gamma_j,\gamma'_j)$.

\section*{2. Lemmas concerning the structure of the root systems
$S$ and $\Phi$.}

\hspace{1.2em} In this section we prove several lemmas. As above,
let $(r_1,\ldots,r_s)$ be the sizes of blocks in the reductive
subalgebra~$\mathfrak{r}$. We denote
$R_k=\displaystyle\sum^{k}_{t=1}r_t$.

\medskip
\textbf{Lemma 2.1.}~\emph{For any root} $(i,j)\in M\setminus S$,
\emph{there exists a number $\widetilde{i}>i$ such that
$(\widetilde{i},j)\in S$ or a number
\raisebox{0pt}[16pt]{$\widetilde{j}<j$} for which}
$(i,\widetilde{j})\in S$.

\medskip
\textbf{Remark.}~In other words, the base $S$ is so structured that
for any cell $(i,j)$ of the matrix in $\mathfrak{m}$, there exists
an element of the base $S$ that lies in the $i$th row and/or in the
$j$th column.

\medskip
\textsc{Proof.}~By definition of the base $S$, for arbitrary
$\gamma=(i,j)$ in $M\setminus S$ there exists a root $\xi=(k,m)\in
S$ such that $\gamma-\xi\in\Delta\!^+$\!, i.e.,
$$\gamma-\xi=\varepsilon_i-\varepsilon_j-(\varepsilon_k-
\varepsilon_m)\in\Delta\!^+\!.$$ The last inclusion holds if $i=k$
or $j=m$.~$\Box$

\medskip
\textbf{Lemma 2.2.}~\emph{Suppose that a root $(i,j)\in S$ is such
that $i\neq R_{k}$ for any} $k$. \emph{Then there exists a number
$\widetilde{j}<j$ such that}
\raisebox{0pt}[16pt]{$(i+1,\widetilde{j})\in S$}.

\medskip
\textsc{Proof.}~Consider the root $(i+1,j)\in M$; it does not lie in
$S$ and beneath it there are no roots from $S$, because otherwise
there are two roots from $S$ in one column. Assume that for any
\raisebox{0pt}[15pt]{$\widetilde{j}<j$} we have
\raisebox{0pt}[13pt]{$(i+1,\widetilde{j})\not\in S$}. Then for the
root $(i+1,j)$ there is no $\xi\in S$ for which
$$(i+1,j)-\xi\in\Delta\!^+\!,$$
which contradicts the definition of the base.~$\Box$

\medskip
\textbf{Corollary.}~\emph{Suppose $(i,j)\in S$ is such that for some
number}~$k$, \emph{we have} $$R_{k-1}<i<R_{k}.$$ \emph{Then for any
$\widetilde{i}$ such that} $i<\widetilde{i}\leqslant R_k$,
\emph{there exists a number $\widetilde{j}<j$ for which}
\raisebox{0pt}[15pt]{$(\widetilde{i},\widetilde{j})\in S$.}

\medskip
\textsc{Proof}~proceeds by induction on the number
$\widetilde{i}$.~$\Box$

\medskip
The following lemma is proved similarly.

\medskip
\textbf{Lemma 2.3.}~\emph{Suppose a root $(i,j)\in S$ is such that
$j\neq R_{k}+1$ for any}~$k$. \emph{Then there exists a number
\raisebox{0pt}[16pt]{$\widetilde{i}>i$} such that the root}
\raisebox{0pt}[17pt]{$(\widetilde{i},j-1)\in S$}.

\medskip
\textbf{Corollary.}~\emph{Assume that $(i,j)\in S$ is such that for
some number}~$k$, \emph{we have}
$$R_{k}+1<j\leqslant R_{k+1}.$$ \emph{Then for any
$\widetilde{j}$ satisfying} $R_k<\widetilde{j}<j$, \emph{there
exists a number $\widetilde{i}>i$ such that}
\raisebox{0pt}[15pt]{$(\widetilde{i},\widetilde{j})\in S$.}

\medskip
\textbf{Lemma 2.4.}~\emph{Let a root $(i,j)\in S$ be such that there
is no root $(\widetilde{i},\widetilde{j})\in S$ with
\raisebox{0pt}[15pt]{$i>\widetilde{i}$} and} $j<\widetilde{j}$.
\emph{Assume that for some $a$ and} $b$, \emph{the following
conditions are satisfied}: $R_a<i\leqslant R_{a+1}$,
$R_{b-1}<j\leqslant R_{b}$. \emph{Then $i=1+R_{a}$ or} $j=R_{b}$.

\medskip
For example, in Diagram 1 the roots $(1,2)$, $(2,10)$, $(9,12)$ and
$(11,14)$ satisfy the assumptions of Lemma 2.4.

\medskip
\textsc{Proof.}~Assume the contrary, i.e., $i>R_{a}+1$ and
$j<R_{b}.$

Consider the root $(i-1,j+1)\in M$; from the conditions of the
lemma, it follows that it does not lie in $S$. Then there exists
$\xi\in S$ for which
$$(i-1,j+1)-\xi\in\Delta\!^+\!.$$
Lemma 2.1 implies that $\xi=(i-1,k)$, where $k<j+1$, or
$\xi=(m,j+1)$, where $i-1<m$. In the first case, by the corollary of
Lemma 2.2, for any number $k$ such that $R_{b-1}<k<j+1$, there
exists a root \mbox{$(i_k,k)\in S$}, where $i_k>i-1$. Then
$\xi\not\in S$, because otherwise in the $k$th column there are two
roots from~$S$, which contradicts the minimality of the base. The
second case is also impossible by the consequence of Lemma 2.3.
Therefore, such a $\xi\in S$ does not exist and $(i-1,j+1)\in S$, a
contradiction.~$\Box$

\medskip
\textbf{Lemma 2.5.}~\emph{Suppose the following $r$ roots lie in}
$S$:
\begin{equation}
(i_1,j-r+1),(i_2,j-r+2),\ldots,(i_r,j)=(i,j),\label{L5_i}
\end{equation}
\emph{where} $i_1>i_2>\ldots>i_r$; \emph{moreover}, \emph{there is
no root $(i_0,j-r)\in S$ such that} $i_0>i_1$. \emph{Then there are
precisely $r$ roots in $S$ of the form}
\begin{equation}
(i+r-1,j_1),(i+r-2,j_2),\ldots,(i,j_r)=(i,j),\label{L5_j}
\end{equation}
\emph{where} $j_1<j_2<\ldots<j_r$.

\emph{The converse is also true}; \emph{namely}, \emph{if $r$ roots
of the form (\ref{L5_j}) lie in $S$ and there is no root
$(i-r,j_0)\in S$ with} $j_0<j_1$, \emph{then there are precisely $r$
roots in $S$ of the form}~(\ref{L5_i}).

\medskip
\textsc{Proof.}~The corollary to Lemma 2.2 asserts that such roots
exist. Let the number of them be equal to $t$:
\begin{equation}
(i+t-1,j_1),(i+t-2,j_2),\ldots,(i,j_t)=(i,j),\quad
j_k<j_{k+1}.\label{L5_j2}
\end{equation}
The first root $(i_1,j-r+1)$ on the list (\ref{L5_i}) is minimal (in
the sense of the relation $\succ$), and the value of $j_1$ is the
least possible. Then $i_1+1=j-r+1,$ consequently, $i_1=j-r$.
Similarly, the first root in (\ref{L5_j2}) is also minimal;
therefore, we have $i+t-1+1=j_1,$ i.e., $j_1=i+t$. Since the root
$\xi=(i,j)\in S$, it follows that a square minor $M_{\xi}$
corresponds to it. We find its size: $i_1-i+1=j-j_1+1,$ whence we
have $j-r-i=j-(i+t)$ and $r=t$.~$\Box$

\medskip
\textbf{Definition 2.6.} Any root $(i,j)\in M$ satisfying the
conditions
$$R_{k-1}<i\leqslant R_k\mbox{ and }R_k<j\leqslant n,$$
is called a \emph{root lying to the right of the $k$th block} in
$\mathfrak{r}$.

\medskip
\textbf{Lemma 2.7.}~\emph{Let roots
$$(i,j_1),(i,j_2),\ldots,(i,j_a)\in S\cup\Phi,$$
where $j_1<j_2<\ldots<j_a$ and} $a>1$, and let other roots from
$S\cup\Phi$ do not lie in the $i$th row. Then
\begin{itemize}
\item[1)] \emph{if} $R_{l-1}<i\leqslant R_{l}$, \emph{then} $j_1=R_l+1$;
\item[2)] \emph{if} $j_1\leqslant R_{k}<j_b$ \emph{for some
number}~$b\leqslant a$, \emph{then} $r_k<b$.
\end{itemize}

\medskip
\textsc{Proof.}
\begin{itemize}
\item[1).] We prove the first item of the lemma. Let
$R_{l-1}<i\leqslant R_{l}$ and $j_1\neq R_l+1$. Since the root
$(i,j_1)\in\Phi$, there exist roots
\raisebox{0pt}[15pt]{$(\widetilde{i},j_1)$}, $\widetilde{i}>i$, and
$(c,i)$ that lie in $S$. Next, since $j_1\neq R_l+1$, by Lemma 2.3
there is a root $(\,\widehat{i},j_1-1)\in S$. Then the pair of roots
\raisebox{0pt}[15pt]{$\big((c,i),(\,\widehat{i},j_1-1)\big)$} is
admissible. Hence we have $(i,j_1-1)\in\Phi$, a contradiction. Thus
$j_1=R_l+1$.

\item[2).] Let we have $R_{m-1}<j_b\leqslant
R_m$ for some number $m$. Show that $r_k<b$ for any $k$ satisfying
the condition $l<k<m$. We prove this by induction on $k$. For the
$(l+1)$th block in $\mathfrak{r}$ the statement is obvious. Suppose
the sizes of blocks in $\mathfrak{r}$ the numbers of which is equal
to $l+1,l+2,\ldots,k-1$ are less than $b$. We show that $r_k<b$.

Suppose the contrary. Let $r_k>b$; then for any number $j$,
where$R_{k-1}<j\leqslant R_k$, a root $(a_j,j)\in S$ exists for some
$a_j\geqslant i$. Thus, there are $r_k$ roots from $S$ of the form
$$(i_1,R_{k-1}+1),(i_2,R_{k-1}+2),\ldots,(i_{r_k},R_k)$$
for some $i_1>i_2>\ldots>i_{r_k}$. Lemma 2.5 implies that there
exist $r_k$ roots in $S$ that lie to the right of a certain block in
$\mathfrak{r}$. Then the size of this block is no less than $r_k$
and thus is greater than $b$. Since its number is less than $k$, we
obtain a contradiction with the inductive assumption. Consequently,
for any number $k$ such that $j_1\leqslant R_{k}<j_b$, we have
$r_k<b$.~$\Box$
\end{itemize}

\section*{3. The root system $T$ and the principal minors}

\hspace{1.2em} Denote by $w_{i,j}:\ \mathfrak{m}\rightarrow K$ the
coordinate function that maps a matrix $x\in\mathfrak{m}$ to the
number that stands at the intersection of the $i$th row and $j$th
column.

Consider $x\in\mathfrak{m}$. Let the corresponding block diagonal
subalgebra $\mathfrak{r}$ consists of blocks of sizes
$(r_1,\ldots,r_s)$. Then $n=r_1+\ldots+r_s.$ We fix the size of the
last block: $r_s=r$. As above, $\displaystyle
R_k=\displaystyle\sum^{k}_{t=1}r_t.$

To prove Theorem 1.7, we need to consider two cases, where there is
a root from $S$ in the $n$th column and there are no roots from $S$
in the $n$th column. To prove the theorem in the first case, we
introduce the notion of principal minor (Definition 3.3) and define
root systems $\Psi$ and $T$ (notation 3.1 and 3.6, respectively).
Throughout the present section, we assume that a parabolic
subalgebra $\mathfrak{p}$ is such that there is a root from $S$ in
the $n$th column. We denote it by $(\widetilde{m},n)$. We take a
number $m\leqslant\widetilde{m}$ such that $(m,n)\in S\cup\Phi$ and
if $i<m$, then $(i,n)\not\in S\cup\Phi$. By Lemma 2.5 and the
corollary to Lemma 2.2, there exists precisely $r$ roots in $S$ of
the form
$$(\widetilde{m}+r-1,l_1),(\widetilde{m}+r-2,l_2),\ldots,
(\widetilde{m}+1,l_{r-1}),(\widetilde{m},l_r)=(\widetilde{m},n)$$
for which $l_1<l_2<\ldots<l_r=n$.

We show graphically the part of the diagram for $\mathfrak{m}$ that
contains these roots:
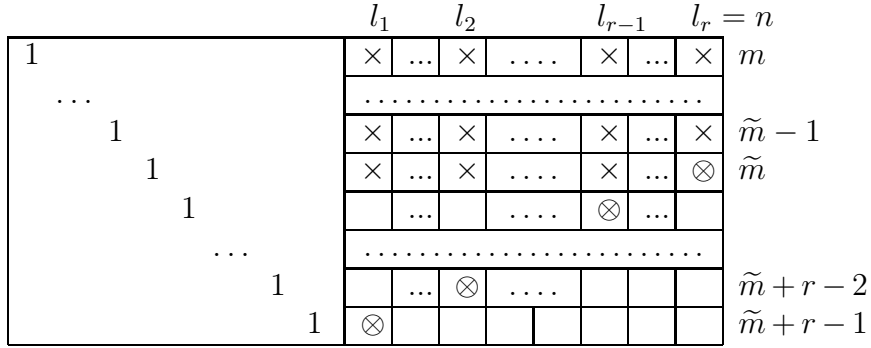
\begin{center}
\begin{tabular}{|p{4cm}|p{0.2cm}|p{0.2cm}|p{0.2cm}|p{0.2cm}|
p{0.2cm}|p{0.2cm}|p{0.2cm}|p{0.2cm}|p{1.7cm}}
\multicolumn{3}{r}{$l_1$\quad\quad}&\multicolumn{3}{l}{$l_2$}&
\multicolumn{2}{l}{$l_{r-1}$}&\multicolumn{2}{l}{$l_r=n$}\\
\cline{1-9} 1&$\times$&...&$\times$&\multicolumn{2}{c|}{$\dotfill$}&
$\times$&...&$\times$&$m$\\
\cline{2-9} \quad$\ldots$&\multicolumn{8}{|c|}{$\dotfill$}&\\
\cline{2-9} \qquad\quad \!\!1&$\times$&...&$\times$&
\multicolumn{2}{c|}{$\dotfill$}&$\times$&...&$\times$&
$\widetilde{m}-1$\\
\cline{2-9} \qquad\qquad \!1&$\times$&...&$\times$&
\multicolumn{2}{c|}{$\dotfill$}&$\times$&...&$\otimes$&
$\widetilde{m}$\\
\cline{2-9} \qquad\qquad\quad 1&&...& &
\multicolumn{2}{c|}{$\dotfill$}&
$\otimes$&...&&\\
\cline{2-9} \qquad\qquad\qquad$\ldots$&
\multicolumn{8}{|c|}{$\dotfill$}&\\
\cline{2-9} \qquad\qquad\qquad\qquad \!1&&...&$\otimes$&
\multicolumn{2}{c|}{$\dotfill$}& & & &$\widetilde{m}+r-2$\\
\cline{2-9} \qquad\qquad\qquad\qquad\quad 1&$\otimes$& & & & & & &
&$\widetilde{m}+r-1$\\
\cline{1-9}
\end{tabular}
\end{center}

Consider all possible chains of roots from $S\cup\Phi$ that contain
the roots $(i,j)$, where $m\leqslant i<m+r$ and
$j\in\{l_1,l_2,\ldots,l_r\}$. Among these chains, we chose those
satisfying the following condition: if $(a,b)$ is a root of such
chain and $R_k<a\leqslant R_{k+1}$ for some $k\geqslant1$, then
\begin{equation}
a\leqslant r+R_k.\label{<=r+R_}
\end{equation}

\textbf{Notation 3.1.}  The set of all roots that occur in such
chains is denoted by~$\Psi$.

\medskip
Obviously, the length of any root chain from $S\cup\Phi$ is less
than the number of blocks in $\mathfrak{r}$. With Example 1, we show
the arrangement of roots from $\Psi$. We mark them by the
symbol~$\boxtimes$.
\begin{center}
{\begin{tabular}{|p{0.1cm}|p{0.1cm}|p{0.1cm}|p{0.1cm}|p{0.1cm}|
p{0.1cm}|p{0.1cm}|p{0.1cm}|p{0.1cm}|p{0.1cm}|p{0.1cm}|p{0.1cm}|
p{0.1cm}|p{0.1cm}|c} \multicolumn{2}{l}{{\small 1\hspace{5pt}
2\!\!}}&\multicolumn{2}{l}{{\small 3\hspace{5pt}
4\!\!}}&\multicolumn{2}{l}{{\small 5\hspace{5pt} 6\!\!}}&
\multicolumn{2}{l}{{\small 7\hspace{5pt}
8\!\!}}&\multicolumn{2}{c}{{\small\makebox[0.2cm][c]{9\ \
10\!}}}&\multicolumn{2}{c}{{\small\makebox[0.2cm][c]{11\hspace{3pt}12}}}&
\multicolumn{2}{c}{{\small\makebox[0.2cm][c]{13\hspace{3pt}14}}}\\
\cline{1-14} 1&$\boxtimes$&&&&&&&&&&&&&{\small 1}\\
\cline{1-14} &\multicolumn{3}{l|}{1}&$\boxtimes$&$\boxtimes$&&&&
$\otimes$&&&&&{\small 2}\\
\cline{5-14} &\multicolumn{3}{c|}{1}&&$\boxtimes$&&&&&&&&
&{\small 3}\\
\cline{5-14} &\multicolumn{3}{r|}{1}&$\otimes$&&&&&&&&&&{\small 4}\\
\cline{2-14} \multicolumn{4}{|c|}{}&\multicolumn{2}{l|}{1}&
$\boxtimes$&&$\boxtimes$&&&&&&{\small 5}\\
\cline{7-14} \multicolumn{4}{|c|}{}&\multicolumn{2}{r|}{1}&
$\boxtimes$&&&&&&&&{\small 6}\\
\cline{5-14} \multicolumn{6}{|c|}{}&1&$\boxtimes$&&&&&&&{\small 7}\\
\cline{7-14} \multicolumn{7}{|c|}{}&\multicolumn{3}{l|}{1}&
$\boxtimes$&$\boxtimes$&&&{\small 8}\\
\cline{11-14} \multicolumn{7}{|c|}{}&\multicolumn{3}{c|}{1}&
$\boxtimes$&$\boxtimes$&&&{\small 9}\\
\cline{11-14} \multicolumn{7}{|c|}{}&\multicolumn{3}{r|}{1}&
$\otimes$&&&&{\small 10}\\
\cline{8-14} \multicolumn{10}{|c|}{}&\multicolumn{2}{l|}{1}&
$\boxtimes$&$\boxtimes$&{\small 11}\\
\cline{13-14} \multicolumn{10}{|c|}{}&\multicolumn{2}{r|}{1}&
$\boxtimes$&&{\small 12}\\
\cline{11-14} \multicolumn{12}{|c|}{}&\multicolumn{2}{l|}{1}&
{\small 13}\\
\multicolumn{12}{|c|}{}&\multicolumn{2}{r|}{1}&{\small 14}\\
\cline{1-14} \multicolumn{14}{c}{Diagram 2}\\
\end{tabular}}
\end{center}

We recall that $s$ is the number of blocks of the reductive
subalgebra $\mathfrak{r}$. The size of the last block in
$\mathfrak{r}$ is equal to $r$. Obviously, for some number
$\widetilde{s}<s$ we have \mbox{$m=R_{\widetilde{s}}+1$.} Lemma 2.7
implies that if $\widetilde{s}+1<t<s$, then the size of the $t$th
block in $\mathfrak{r}$ is less than the size of the last block $r$.
For this reason, all roots $(i,j)$ lying in $S\cup\Phi$ such that
$i>R_{\widetilde{s}+1}$ are contained in the root system $\Psi$.

\medskip
\textbf{Example 2.} We give one more example of the diagram of a
parabolic subalgebra with sizes of diagonal blocks
$(2,2,1,3,2,1,3)$. Here $\widetilde{s}=3$ and $m=R_{3}+1=6$.
\begin{center}
{\begin{tabular}{|p{0.1cm}|p{0.1cm}|p{0.1cm}|p{0.1cm}|p{0.1cm}|
p{0.1cm}|p{0.1cm}|p{0.1cm}|p{0.1cm}|p{0.1cm}|p{0.1cm}|p{0.1cm}|
p{0.1cm}|p{0.1cm}|c} \multicolumn{2}{l}{{\small 1\hspace{5pt}
2\!}}&\multicolumn{2}{l}{{\small 3\hspace{5pt}
4\!\!}}&\multicolumn{2}{l}{{\small 5\hspace{5pt} 6\!\!}}&
\multicolumn{2}{l}{{\small 7\hspace{5pt}
8\!\!}}&\multicolumn{2}{c}{{\small\makebox[0.2cm][c]{9\ \
10\!}}}&\multicolumn{2}{c}{{\small\makebox[0.2cm][c]{11\hspace{3pt}12}}}&
\multicolumn{2}{c}{{\small\makebox[0.2cm][c]{13\hspace{3pt}14}}}\\
\cline{1-14} \multicolumn{2}{|l|}{1}&&$\boxtimes$&&&&&&&&&&
&{\small 1}\\
\cline{3-14}\multicolumn{2}{|r|}{1}&$\boxtimes$&&&&&&&&&&&
&{\small 2}\\
\cline{1-14}\multicolumn{2}{|c|}{}&\multicolumn{2}{|l|}{1}&
$\boxtimes$&&$\boxtimes$&&&&&&&&{\small 3}\\
\cline{5-14} \multicolumn{2}{|c|}{}&\multicolumn{2}{r|}{1}&
$\boxtimes$&&&&&&&&&&{\small 4}\\
\cline{3-14} \multicolumn{4}{|c|}{}&1&$\boxtimes$&&&&&&&&&
{\small 5}\\
\cline{5-14} \multicolumn{5}{|c|}{}&\multicolumn{3}{l|}{1}&
$\boxtimes$&$\boxtimes$&&&&$\boxtimes$&{\small 6}\\
\cline{9-14} \multicolumn{5}{|c|}{}&\multicolumn{3}{c|}{1}&
$\boxtimes$&$\boxtimes$&&&&&{\small 7}\\
\cline{9-14} \multicolumn{5}{|c|}{}&\multicolumn{3}{r|}{1}&
$\boxtimes$&&&&&&{\small 8}\\
\cline{6-14} \multicolumn{8}{|c|}{}&\multicolumn{2}{l|}{1}&
$\boxtimes$&&$\boxtimes$&&{\small 9}\\
\cline{11-14} \multicolumn{8}{|c|}{}&\multicolumn{2}{r|}{1}&
$\boxtimes$&&&&{\small 10}\\
\cline{9-14} \multicolumn{10}{|c|}{}&1&$\boxtimes$&&&
{\small 11}\\
\cline{11-14} \multicolumn{11}{|c|}{}&\multicolumn{3}{l|}{1}&
{\small 12}\\
\multicolumn{11}{|c|}{}&\multicolumn{3}{c|}{1}&{\small 13}\\
\multicolumn{11}{|c|}{}&\multicolumn{3}{r|}{1}&{\small 14}\\
\cline{1-14} \multicolumn{14}{c}{Diagram 3}\\
\end{tabular}}
\end{center}

Denote $\widetilde{r}_k=\min(r,r_k)$. Then from the definition of
the set $\Psi$ it follows that if $(i,j)\in\Psi$, then for $k>1$ we
have
$$R_{k-1}<i\leqslant\widetilde{r}_k+R_{k-1}.$$

\medskip
\textbf{Lemma 3.2.} \emph{Let} $k>1$. \emph{The following statement
are valid}:
\begin{enumerate}
\item \emph{There are $\widetilde{r}_k$ roots from $S\cup\Phi$ in the row
with number} \mbox{$R_{k-1}+1$}.
\item \emph{The first} (\emph{according to the number of a column})
$\widetilde{r}_k$ \emph{roots from $S\cup\Phi$ in the row with
number \mbox{$R_{k-1}+1$} lie in} $\Psi$.
\item \emph{All rots from $\Psi$ that lie to the right of the $k$th block
in $\mathfrak{r}$ are contained in the rows with numbers
$R_{k-1}+1,R_{k-1}+2,\ldots,R_{k-1}+\widetilde{r}_k$ and in the same
columns as the roots from $\Psi$ in the row with number}
$R_{k-1}+1$.
\end{enumerate}

\medskip
\textsc{Proof.}
\begin{enumerate}
\item
\begin{enumerate}
\item First we prove that at least $\widetilde{r}_k$ roots from $S$
lie to the right of the $k$th block.

The proof proceeds by contradiction. assume that to the right of the
$k$th block, there are less than $\widetilde{r}_k$ roots from $S$.
Then there are no roots from $S$ in the row with the number
$R_{k-1}+1$. Otherwise to the right of the $k$th block there are
$r_k\geqslant\widetilde{r}_k$ roots from $S$. By Lemma 2.1, for each
root $(R_{k-1}+1,j)\in M$ lying in the row with number $R_{k-1}+1$,
a root $\xi=(i,j)\in S$ with $i>R_{k-1}+1$ exists.

Obviously, $r=r_s\geqslant\widetilde{r}_k$. Let $l$ be the least
number for which $r_l\geqslant\widetilde{r}_k$, $l>k$. Since
$r_s=r\geqslant\widetilde{r}_k$ and $s>k$, we have $l\leqslant s$.

Next, in each of the columns with numbers
$$R_{l-1}+1,R_{l-1}+2,\ldots,R_{l-1}+r_l=R_l$$
there is a root $(a_i,R_{l-1}+i)\in S$, where $i=1,\ldots,r_l$.
Lemma 2.3 implies that $a_1>a_2>\ldots>a_{r_l}$. Then, by Lemma 2.5,
there exists a number $k\leqslant p<l$ such that there are $r_l$
roots from $S$ to the right of the $p$th block in $\mathfrak{r}$.
Then $r_p\geqslant r_l\geqslant\widetilde{r}_k$. Since $l$ was
chosen to be minimal, we obtain $l=k$. But this contradicts the
assumption that to the right of the $k$th block there are less than
$\widetilde{r}_k$ roots from $S$. Thus, to the right of the $k$th
block there are $\widetilde{r}_k$ roots from $S$. Part (a) is
proved.

\item We complete the proof of statement 1. Let the size of the
$k$th block in $\mathfrak{r}$ be greater than 1. For $r_k=1$ the
proof is trivial, because $(R_k,R_{k}+1)\in S$ for any number $k$.

The corollary to Lemma 2.3 states that the roots from $S$ lying to
the right of the $k$th block in $\mathfrak{r}$ have the form
$$\xi_{1}=(R_{k},j_1),\xi_{2}=
(R_{k}-1,j_2),\ldots,\xi_{\widetilde{r}_k}=
(R_{k}-\widetilde{r}_k+1,j_{\widetilde{r}_k}),\ldots,$$ where
$j_1=R_k+1$ and $j_1<j_2<\ldots<j_{\widetilde{r}_k}$. We show that
all the roots of the form $(R_{k-1}+1,j_a)$, where $1\leqslant
a\leqslant\widetilde{r}_k$, lie in $S\cup\Phi$.

Consider the root $\eta=(R_{k-1},R_{k-1}+1)$; since $k>1$, this root
obviously lies in $S$. Denote
$$I=\left\{\begin{array}{l}
1,2,\ldots,\widetilde{r}_k-1\mbox{, если }\widetilde{r}_k=r_k;\\
1,2,\ldots,\widetilde{r}_k\mbox{, если }\widetilde{r}_k<r_k.\\
\end{array}\right.$$
Every root $\gamma_a=(R_{k-1}+1,R_{k}-a+1)$, where $a\in I$, lies
in~$\Delta^{\!+}_{\mathfrak{r}}$. Consequently, all pairs of roots
$(\eta,\xi_a)$, $a\in I$, are admissible. Denote
\begin{equation}
\zeta_a=\left\{\begin{array}{l}
\gamma_a+\xi_a=(R_{k-1}+1,j_a)\mbox{, если
}\widetilde{r}_k<r_k,\\
\zeta_{\widetilde{r}_k}=\xi_{\widetilde{r}_k}\mbox{, если
}\widetilde{r}_k=r_k.\\
\end{array}\right.\label{zeta}
\end{equation}
Then the roots $\zeta_a$ are in $\Phi$ for any $a\in I$, which
proves the first statement of the lemma.
\end{enumerate}
\end{enumerate}

We prove statement 2 and 3 for numbers $k$ that do not exceed
$\widetilde{s}+1$. For remaining $k$, there is nothing to prove.

\begin{enumerate}
\item[2.]  The proof of the second statement proceeds by induction
on the number $k$, starting with the greatest one. For
$k=\widetilde{s}+1$ the statement is obvious. Suppose statement 2 is
valid for roots from $\Psi$ that lie to the right of the blocks in
$\mathfrak{r}$ with numbers $k+1,k+2,\ldots,\widetilde{s}+1$. we
show this for the block with number $k$; namely, we prove that for
any $1\leqslant a\leqslant\widetilde{r}_k$, the root $\zeta_a$,
defined in~(\ref{zeta}), lie in the root system~$\Psi$.

Lemma 2.5 implies that there exist roots in $S$ of the form
$$(i_1,j_a-a+1),(i_2,j_a-a+2),\ldots,(i_a,j_a)=(R_{k}-a+1,j_a)$$
for some $i_1>i_2>\ldots>i_a=R_{k}-a+1$. For the first root on the
list, we have $i_1=j_a-a=R_b$ for some number $b\geqslant k$.

Denote $\xi=(R_{b+1},R_{b+1}+1)$. It is easy to see that $\xi\in S$.
If $j_a<R_{b+1}$, then
$(j_a,R_{b+1})\in\Delta^{\!+}_{\mathfrak{r}}$; then the pair
$(\xi_a,\xi)$ is admissible. Consequently, the root
$(j_a,R_{b+1}+1)$ lies in $\Phi$. If $j_a=R_{b+1}$, then
$(j_a,R_{b+1}+1)$lies in the base $S$, i.e., in either case
($j_a<R_{b+1}$ and~$j_a=R_{b+1}$), the roots $(j_a,R_{b+1}+1)\in
S\cup\Phi$. It remains to prove that $(j_a,R_{b+1}+1)\in\Psi$.
Consider the chain
$$(j_a,R_{b+1}+1),(R_{b+1}+1,R_{b+2}+1),(R_{b+2}+1,R_{b+3}+1),\ldots,
(R_{\widetilde{s}}+1,R_{\widetilde{s}+1}+1),$$ all of its roots lie
in $S\cup\Phi$. Since
$$j_a=a+R_b\leqslant\widetilde{r}_k+R_b\leqslant r+R_b,$$
the pair $(j_a,R_{b+1}+1)$satisfies (\ref{<=r+R_}). Consequently
$(j_a,R_{b+1}+1)\in\Psi$. Hence the root $\zeta_a$ also lies in
$\Psi$.

\item[3.] We prove statement 3. Consider a root $(i,j)$, where
$R_{k-1}+1<i\leqslant R_{k-1}+\widetilde{r}_k$, $j_{a-1}<j<j_a$ for
some $a\in\{2,\ldots,\widetilde{r}_k\}$. Suppose $(i,j)\in
S\cup\Psi$. Therefore, there is a root $(\widetilde{i},j)\in S$,
$\widetilde{i}>i$. Then, similarly to item 1 and 2 of the lemma, we
obtain $(R_{k-1}+1,j)\in\Psi$, a contradiction with $j\neq j_a$ for
any $a$. Thus, $(i,j)\not\in\Psi$.

It remains to show that there are no roots from $\Psi$ that lie in
the rectangle $R_{k-1}<i\leqslant R_{k-1}+\widetilde{r}_k$,
$j_{\widetilde{r}_k}<j\leqslant n$. Assume the contrary. Let
$(i,j)\in\Psi$. Then $(i,j)\in S\cup\Phi$ and there exists a number
$a>\widetilde{r}_k$ such that $\xi_a=(R_{k}-a+1,j)\in S$, where
$R_k-a+1\geqslant i$.

Note that if $a>\widetilde{r}_k$, then $r_k>\widetilde{r}_k$, whence
$\widetilde{r}_k=r$ and $a>r$. Repeating the arguments used in item
2, we obtain the following estimate for $j$:
$$j=R_b+a>r+R_b.$$
Consequently, there are no roots from $\Psi$ in the $j$th row. Then
there are no roots from $\Psi$ in the $j$th column as well, a
contradiction with $(i,j)\in\Psi$.~$\Box$
\end{enumerate}

To each block, except for the last one, of the reductive part
$\mathfrak{r}$ we put in correspondence a certain minor of the
matrix $\mathbb{X}$ in accordance with the following rule. Let
$j_1,j_2,\ldots,j_{\widetilde{r}_k}$ be the numbers of columns in
which there are all roots from $\Psi$ that lie to the right of the
$k$th block in~$\mathfrak{r}$.

\medskip
\vbox{\textbf{Definition 3.3.} The minor in $\mathbb{X}$ that stands
at the intersection of $R_{k-1}<i\leqslant\widetilde{r}_k+R_{k-1}$
rows and $j_1,j_2,\ldots,j_{\widetilde{r}_k}$ columns is called a
\emph{principal minor}.

\medskip
We indicate all principal minors in Example 1 and 2. As above, the
minor $M_I^J$ is a minor of the matrix $\mathbb{X}$ with the ordered
systems of rows $I$ and columns $J$.

In Example 1, the minors  $M_1^2$, $M_{2,3}^{5,6}$, $M_{5,6}^{7,9}$,
$M_7^8$, $M_{8,9}^{11,12}$, and
\raisebox{0pt}[17pt]{$M_{11,12}^{13,14}$} are principal. In Example
2, $M_{1,2}^{3,4}$, $M_{3,4}^{5,7}$, $M_{5}^{6}$,
$M_{6,7,8}^{9,10,14}$, $M_{9,10}^{11,12}$, and
\raisebox{0pt}[17pt]{$M_{11}^{12}$} are principal.

\medskip
The size of the principal minor that stands to the right of the
$k$th block in the subalgebra~$\mathfrak{r}$ is equal to
$\widetilde{r}_k=\mathrm{min}(r,r_k)$. The number of the principal
minors is one less than the number of blocks in $\mathfrak{r}$. Note
that no two principal minors have a common column/row.

\medskip
\textbf{Lemma 3.4.} \emph{Let $(i,j)$ be a root from} $\Psi$.
\emph{Then for any $\widetilde{j}<j$ such that}
$(i,\widetilde{j})\not\in\Delta^{\!+}_{\mathfrak{r}}$, \emph{a root}
$(\widetilde{i},\widetilde{j})$}, $\widetilde{i}\geqslant i$,
\emph{lying in $\Psi$ exists}.

\medskip
\textsc{Proof.} Let
$(i,\widetilde{j})\not\in\Delta_{\mathfrak{r}}^{\!+}$ and
$\widetilde{j}<j$. By definition of the base $S$, there is a root
$\xi\in S$ such that
\raisebox{0pt}[15pt]{$(i,\widetilde{j})-\xi\in\Delta^{\!+}\!$}.
Lemma 2.1 implies that $\xi=(i,b)$, where
\raisebox{0pt}[12pt]{$b<\widetilde{j}$}, or $\xi=(a,\widetilde{j})$,
where $a>i$. The first case os impossible, because in a row, there
is no root from $\Phi$ to the right of a root from $S$; therefore
$\xi=(a,\widetilde{j})$, where $a>i$.

For some number $l$, we have $R_{l-1}<a\leqslant R_l$. The first
assertion of Lemma 3.2 states that
\raisebox{0pt}[15pt]{$(R_{l-1}+1,\widetilde{j})\in S\cup\Phi$.} we
show that this root lies in the system $\Psi$.

For some $k$ we have $R_{k-1}<\widetilde{j}\leqslant R_k$.
Obviously, $(R_k,R_k+1)\in S$. If
\raisebox{0pt}[15pt]{$\widetilde{j}<R_k$,} then
\raisebox{0pt}[12pt]{$(\widetilde{j},R_k)\in\Delta_{\mathfrak{r}}^{\!+}$}
and the pair of roots $\big((a,\widetilde{j}),(R_k,R_k+1)\big)$ is
admissible. Consequently,
\raisebox{0pt}[12pt]{$(\widetilde{j},R_k)\in\Phi$}. If
$\widetilde{j}=R_k$, then $(\widetilde{j},R_k+1)\in S$. Thus,
\raisebox{0pt}[17pt]{$(\widetilde{j},R_k+1)\in S\cup\Phi$} for any
\raisebox{0pt}[15pt]{$\widetilde{j}$} such
that($R_{k-1}<\widetilde{j}\leqslant R_k$. Lemma 3.2
and(\ref{<=r+R_}) imply that all the roots in the chain
$$(\widetilde{j},R_{k}+1),(R_k+1,R_{k+1}+1),\ldots,
(R_{\widetilde{s}}+1,R_{\widetilde{s}+1}+1)$$ lie in $\Psi$. Next,
from Lemma 2.7 it follows that the size of the $k$th block in
$\mathfrak{r}$ is less than~$r$. Then
$$\widetilde{j}\leqslant R_k=R_{k-1}+r_k<R_{k-1}+r.$$
From (\ref{<=r+R_}) it follows that the root
$(R_{l-1}+1,\widetilde{j})$ is also contained in~$\Psi$.

Thus, we showed that for any $\widetilde{j}<j$ there is a number
$\widetilde{i}=R_{l-1}+1$ such that
\raisebox{0pt}[15pt]{$(\widetilde{i},\widetilde{j})\in\Psi$}.~$\Box$

\medskip
\textbf{Notation 3.5.} Denote by $\mathcal{X}$ the set of matrices
$x$ in $\mathfrak{m}$ satisfying the following conditions:
\begin{enumerate}
\item The principal minors are triangular,
$$\left|\begin{array}{cccc}
*&*&\ldots&*\\
*&*&\ldots&0\\
\ldots&\ldots&\ldots&\ldots\\
*&0&\ldots&0\\
\end{array}\right|,$$
and the numbers on the secondary diagonal are different from zero.
In the above diagram we denote by asterisks any numbers from the
ground field $K$.

\item If $(i,j)\in\Psi$ and $j$ is the greatest number for a given
$i$, then $w_{i,\widetilde{j}}(x)=0$ for any $\widetilde{j}>j$.

\item Let a principal minor stand at the intersection of the rows
$i,i+1,\ldots,i+\widetilde{r}-1$ and columns
$j_1<j_2<\ldots<j_{\widetilde{r}}$. Let
$$(\widetilde{i},\widetilde{j})\not\in\Psi\mbox{, where
}i\leqslant\widetilde{i}<i+\widetilde{r}\mbox{ and
}j_{a-1}<\widetilde{j}<j_a$$ for some
$a\in\{2,\ldots\widetilde{r}\}$. Then if the root $(i,j_a)$ stands
in the principal minor below the secondary diagonal (in this case
$w_{\widetilde{i},j_a}(x)=0$), then
$w_{\widetilde{i},\widetilde{j}}(x)=0$.
\end{enumerate}

We present an example of the set $\mathcal{X}\subset\mathfrak{m}$ in
which diagram 1 is a diagram for~$\mathfrak{m}$.

$$\left(\begin{tabular}{p{0.1cm}|p{0.1cm}|p{0.1cm}|p{0.1cm}|
p{0.1cm}|p{0.1cm}|p{0.1cm}|p{0.1cm}|p{0.1cm}|p{0.1cm}|p{0.1cm}
|p{0.1cm}|p{0.1cm}|p{0.1cm}}
0&\multicolumn{2}{l}{*\quad0}&\multicolumn{2}{l}{0\quad0}&
\multicolumn{2}{l}{0\quad0}&\multicolumn{2}{l}{0\quad0}&
\multicolumn{2}{l}{0\quad0}&\multicolumn{3}{l}{0\quad0\quad0}\\
\cline{1-4}
&\multicolumn{3}{l|}{0\quad0\quad0}&\multicolumn{2}{l}{*\quad*}&
\multicolumn{2}{l}{0\quad0}&\multicolumn{2}{l}{0\quad0}&
\multicolumn{2}{l}{0\quad0}&\multicolumn{2}{l}{0\quad0}\\
&\multicolumn{3}{l|}{0\quad0\quad0}&\multicolumn{2}{l}{*\quad0}&
\multicolumn{2}{l}{0\quad0}&\multicolumn{2}{l}{0\quad0}&
\multicolumn{2}{l}{0\quad0}&\multicolumn{2}{l}{0\quad0}\\
&\multicolumn{3}{l|}{0\quad0\quad0}&\multicolumn{2}{l}{*\quad*}&
\multicolumn{2}{l}{*\quad*}&\multicolumn{2}{l}{*\quad*}&
\multicolumn{2}{l}{*\quad*}&\multicolumn{2}{l}{*\quad*}\\
\cline{2-6} \multicolumn{4}{l}{}&\multicolumn{2}{|l|}{0\quad0}&
\multicolumn{2}{l}{*\quad*}&\multicolumn{2}{l}{*\quad0}&
\multicolumn{2}{l}{0\quad0}&\multicolumn{2}{l}{0\quad0}\\
\multicolumn{4}{l}{}&\multicolumn{2}{|l|}{0\quad0}&
\multicolumn{2}{l}{*\quad0}&\multicolumn{2}{l}{0\quad0}
&\multicolumn{2}{l}{0\quad0}&\multicolumn{2}{l}{0\quad0}\\
\cline{5-7} \multicolumn{6}{l|}{}&0&
\multicolumn{3}{l}{*\quad0\quad0}&\multicolumn{2}{l}{0\quad0}
&\multicolumn{2}{l}{0\quad0}\\
\cline{7-10}
\multicolumn{7}{l}{}&\multicolumn{3}{|l|}{0\quad0\quad0}&
\multicolumn{2}{l}{*\quad*}&\multicolumn{2}{l}{0\quad0}\\
\multicolumn{7}{l}{}&\multicolumn{3}{|l|}{0\quad0\quad0}&
\multicolumn{2}{l}{*\quad0}&\multicolumn{2}{l}{0\quad0}\\
\multicolumn{7}{l}{}&\multicolumn{3}{|l|}{0\quad0\quad0}&
\multicolumn{2}{l}{*\quad*}&\multicolumn{2}{l}{*\quad*}\\
\cline{8-12} \multicolumn{10}{c}{{\Large0}\qquad\qquad\qquad\qquad}&
\multicolumn{2}{|l|}{0\quad0}&
\multicolumn{2}{l}{*\quad*}\\
\multicolumn{10}{c}{}&\multicolumn{2}{|l|}{0\quad0}&
\multicolumn{2}{l}{*\quad0}\\
\cline{11-14}
\multicolumn{12}{l}{}&\multicolumn{2}{|l}{0\quad0}\\
\multicolumn{12}{l}{}&\multicolumn{2}{|l}{0\quad0}
\end{tabular}\right)$$

\medskip
\textbf{Notation 3.6.} Denote by $T$ a set of roots
$(i,j)\in\Delta\!^+$ that satisfy the following conditions:
\begin{enumerate}
\item There exists a number $\widetilde{i}>i$ such that
$(\widetilde{i},j)\in\Psi$.
\item In the $i$th row there are no roots from $\Psi$.
\end{enumerate}
Denote $\mathcal{T}=\mathrm{Span}\left\{E_{\xi},\ \xi\in T\right\}$.
Below we give the diagram of Example 1; here the roots from $T$ are
denoted y the symbol $\bullet$.
\begin{center}
{\begin{tabular}{|p{0.1cm}|p{0.1cm}|p{0.1cm}|p{0.1cm}|p{0.1cm}|
p{0.1cm}|p{0.1cm}|p{0.1cm}|p{0.1cm}|p{0.1cm}|p{0.1cm}|p{0.1cm}|
p{0.1cm}|p{0.1cm}|c} \multicolumn{2}{l}{{\small 1\hspace{5pt}
2\!\!}}&\multicolumn{2}{l}{{\small 3\hspace{5pt}
4\!\!}}&\multicolumn{2}{l}{{\small 5\hspace{5pt} 6\!\!}}&
\multicolumn{2}{l}{{\small 7\hspace{5pt}
8\!\!}}&\multicolumn{2}{c}{{\small\makebox[0.2cm][c]{9\ \
10\!}}}&\multicolumn{2}{c}{{\small\makebox[0.2cm][c]{11\hspace{3pt}12}}}&
\multicolumn{2}{c}{{\small\makebox[0.2cm][c]{13\hspace{3pt}14}}}\\
\cline{1-14} 1&$\boxtimes$&&&&&&&&&&&&&{\small 1}\\
\cline{1-14} &\multicolumn{3}{l|}{1}&$\boxtimes$&$\boxtimes$&&&&
$\otimes$&&&&&{\small 2}\\
\cline{5-14} &\multicolumn{3}{c|}{1}&&$\boxtimes$&&&&&&&&&
{\small 3}\\
\cline{5-14} &\multicolumn{3}{r|}{1}&$\otimes$&&$\bullet$&
$\bullet$&$\bullet$&&$\bullet$&$\bullet$&$\bullet$&$\bullet$&
{\small 4}\\
\cline{2-14} \multicolumn{4}{|c|}{}&\multicolumn{2}{l|}{1}&
$\boxtimes$&&$\boxtimes$&&&&&&{\small 5}\\
\cline{7-14} \multicolumn{4}{|c|}{}&\multicolumn{2}{r|}{1}&
$\boxtimes$&&&&&&&&{\small 6}\\
\cline{5-14} \multicolumn{6}{|c|}{}&1&$\boxtimes$&&&&&&&
{\small 7}\\
\cline{7-14} \multicolumn{7}{|c|}{}&\multicolumn{3}{l|}{1}&
$\boxtimes$&$\boxtimes$&&&{\small 8}\\
\cline{11-14} \multicolumn{7}{|c|}{}&\multicolumn{3}{c|}{1}&
$\boxtimes$&$\boxtimes$&&&{\small 9}\\
\cline{11-14} \multicolumn{7}{|c|}{}&\multicolumn{3}{r|}{1}&
$\otimes$&&$\bullet$&$\bullet$&{\small 10}\\
\cline{8-14} \multicolumn{10}{|c|}{}&\multicolumn{2}{l|}{1}&
$\boxtimes$&$\boxtimes$&{\small 11}\\
\cline{13-14} \multicolumn{10}{|c|}{}&\multicolumn{2}{r|}{1}&
$\boxtimes$&&{\small 12}\\
\cline{11-14} \multicolumn{12}{|c|}{}&\multicolumn{2}{l|}{1}&
{\small 13}\\
\multicolumn{12}{|c|}{}&\multicolumn{2}{r|}{1}&{\small 14}\\
\cline{1-14} \multicolumn{14}{c}{Diagram 4}\\
\end{tabular}}
\end{center}

\section*{4. Canonical matrices on $N$-orbits of general position}

\hspace{1.2em} First we state two propositions, which will be
necessary for the proof of Theorem 1.7.

\medskip
\textbf{Remark 4.1.} The adjoint action
$$\mathrm{Ad}_{g(t)}\,x=g(t)\cdot x\cdot g(t)^{-1}$$
on an arbitrary matrix $x$ by the one-parameter subgroup
$$g(t)=E+t\cdot E_{\psi}\subset N,\ \psi=(u,v)\in\Delta\!^+$$
reduces to the composition of two transformations:
\begin{itemize}
\item[1)] the row with number $v$ multiplied by~$t$ is added to the
row $u$ of the matrix $x$;
\item[2)]  the column with the number $u$ multiplied by $-t$ is
added to the column with the number $v$ of the matrix $x$.
\end{itemize}

\medskip
\textbf{Proposition 4.2.} \emph{For any $x\in\mathcal{X}$ and any
matrix} $y\in\mathcal{T}$, \emph{there is $g\in N$ such that}
$\mathrm{Ad}_{g}(x+y)=x$.

\medskip
\textsc{Proof.} It is sufficient to prove the following assertion.
Let $j$ be the number of a column such that all entries of the
matrix $y$ in the columns with numbers greater that $j$ are zero. We
denote by $\{i_1,\ldots,i_p\}$ the numbers of the rows for which the
entries of the matrix $y$ in the $j$th column are equal to zero. We
claim that there exists $g\in N$ such that
$$\mathrm{Ad}_g(x+y)=x+\widetilde{y},$$
where $\widetilde{y}\in\mathcal{T}$ and all the entries of
$\widetilde{y}$ that lie in the columns with numbers at least $j$
are zero.

Let $i\in\{i_1,\ldots,i_p\}$; since $(i,j)\in T$, a root
$(\widetilde{i},j)\in\Psi$, $\widetilde{i}>i$, exists. Then there
exists a principal minor $M_I^J$, where $I$ is the set of rows and
$J$ is the set of columns of the minor, such that
\raisebox{0pt}[12pt]{$\widetilde{i}\in I$} and $j\in J$. Let the
principal minor $M_I^J$ stand to the right of the $(c+1)$th block in
$\mathfrak{r}$; denote $k=R_c+1$ and, as above,
$\widetilde{r}=\widetilde{r}_{c+1}=\min(r,r_{c+1})$. Let
\begin{center}
$I=\{k,k+1,\ldots,k+\widetilde{r}-1\}$ and
$J=\{j_1,j_2,\ldots,j_{\widetilde{r}}\}$
\end{center}
for some $j_1<j_2<\ldots<j_{\widetilde{r}}$. Denote by $l$ the
number for which $j=j_l$. Note that the root
$(k+\widetilde{r}-l,j_l)$ stands on the secondary diagonal of a
principal minor, a nd thus the value of $w_{k+\widetilde{r}-l,j}(x)$
does not equal to zero.

Let $i\in\{i_1,\ldots,i_p\}$; then $w_{i,j}(y)\neq0$. Consider the
adjoint action on $x+y$ of the element
$$g=\prod_{i=i_1,\ldots,i_a}g_i(t_i),$$
where $g_i(t_i)=E+t_i\cdot E_{i,k+\widetilde{r}-l}$ and
$t_i=-\displaystyle\frac{w_{i,j}(y)}{w_{k+\widetilde{r}-l,j}(x)}$.

Then the columns in the matrices $\mathrm{Ad}_g(x+y)$ and $x$ the
numbers of which are at least $j$ coincide. Therefore, the entries
of the matrix
$$\widetilde{y}=\mathrm{Ad}_g(x+y)-x,$$
that lie in the columns with numbers at least $j$ are zero. Hence we
obtain $\mathrm{Ad}_g(x+y)=x+\widetilde{y}$.

To complete the proof, it remains to prove that
$\widetilde{y}\in\mathcal{T}$. We show that if
$w_{a,b}(\widetilde{y})\neq0$, then $(a,b)\in T$.

Thus, let $w_{a,b}(\widetilde{y})\neq0$ and
$w_{a,b}(\widetilde{y})\neq w_{a,b}(y)$ for a certain root $(a,b)$.
By Remark 4.1, either $a=i$ or $b=k+\widetilde{r}-l$.

\begin{enumerate}
\item Let $w_{i,b}(\widetilde{y})\neq0$ and $w_{i,b}(\widetilde{y})\neq
w_{i,b}(y)$, where $b\in\{j_1,j_1+1,\ldots,j\}$. We show that
$(i,b)\in T$. Lemma 3.4 implies that for every number $b<j$, a root
$(i_b,b)\in\Psi$with number $i_b>i$ exists. This, together with the
fact that there are no roots from $\Psi$ in the $i$th row, implies
that the root $(i,b)$ lies in $T$.

\item Let $w_{a,k+\widetilde{r}-l}(\widetilde{y})\neq0$ and
$w_{a,k+\widetilde{r}-l}(\widetilde{y})\neq
w_{a,k+\widetilde{r}-l}(y)$ for some $a$. we show that
$(a,k+\widetilde{r}-l)\in T$.
\begin{enumerate}
\item We prove that there are no roots from $\Psi$ in the row with
number $a$. To this end, we first prove that there are no roots from
$\Psi$ in the $i$th column. Assume the contrary. Then the $i$th row
would also contain roots from $\Psi$, which contradicts the
definition of the root system~$T$. Consequently, in the column with
number $i$ there are no roots from $\Psi$. Next, since
$w_{a,k+\widetilde{r}-l}(y)\neq
w_{a,k+\widetilde{r}-l}(\widetilde{y})$, it follows that the entry
$(a,k+\widetilde{r}-l)$ of the matrix $x+y$ is changed under the
action of $g$, and thus $w_{a,i}(x+y)\neq0$. Now assume that in the
$a$th row there is a root $(a,v)$ from $\Psi$, and let $v$ be the
greatest possible. Then, since in the matrices from $\mathcal{X}$
and $\mathcal{T}$ there are zeros to the right of the principal
minors, we obtain $i\leqslant v$. by Lemma 3.4, there is a root from
$\Psi$ in the $i$th column, a contradiction. Therefore, there are no
roots from $\Psi$ in the $a$th row.

\item Now we show that there is a number $\widetilde{a}>a$ such that
$(\widetilde{a},k+\widetilde{r}-l)\in\Psi$. Assume that there are no
roots from $\Psi$ in the $(k+\widetilde{r}-l)$th column. Then there
are no roots from $S$ in the $(k+\widetilde{r}-l)$th column.
Otherwise Lemma 3.2 and the fact that the row with number
$(k+\widetilde{r}-l)$ contains roots from the system $\Psi$ imply
that there is a root $(\widetilde{a},k+\widetilde{r}-l)\in\Psi$ for
some $\widetilde{a}$. But in the $(k+\widetilde{r}-l)$th column
there are no roots from $S$ only if for every
$\widetilde{c}\leqslant c$ the size of the $\widetilde{c}$th block
in $\mathfrak{r}$ is less than the size of the last block. In
particular, the size of the block to the right of which the root
$(a,k+\widetilde{r}-l)$ stands is less than~$r$. Then the row with
number $a$ is a row of a certain principal minor, and since
$w_{a,i}(x+y)\neq0$, there is a root from $\Psi$ in the $a$th row,
which contradicts item (a). Therefore, the $(k+\widetilde{r}-l)$th
column contains a root $(\widetilde{a},k+\widetilde{r}-l)\in\Psi$.

Now assume that $\widetilde{a}<a$. The root
$(\widetilde{a},k+\widetilde{r}-l)$ lies in $\Psi$ and
$i<k+\widetilde{r}-l$. Since
$(a,i)\not\in\Delta_{\mathfrak{r}}^{\!+}$, we have
$(\widetilde{a},i)\not\in\Delta_{\mathfrak{r}}^{\!+}$. Then the
conditions of Lemma 3.4 are fulfilled, whence it follows that there
is a root from $\Psi$ in the $i$th column, but this contradicts
item~(a). Consequently, $\widetilde{a}>a$.
\end{enumerate}

Items (a) and (b) imply that the root $(a,k+\widetilde{r}-l)\in T$.
\end{enumerate}

Thus, all the roots of the form $(a,b)$, where
$w_{a,b}(\widetilde{y})\neq0$, lie in the root system~$T$.~$\Box$

\medskip
\textbf{Proposition 4.3.} \emph{Assume that a parabolic subalgebra
$\mathfrak{p}$ is such that there is a root from $S$ in the last
column}. \emph{Then there exists a nonempty Zariski-open subset $U$
in $\mathcal{Y}$ such that the $N$-orbit of any matrix from $U$ has
a nonzero intersection with} $\mathcal{X}$.

\medskip
\textsc{Proof.} Consider the linear span
\begin{equation}
\mathcal{Y}+\sum K\cdot E_{\eta},\label{sum}
\end{equation}
where summation proceeds over all roots $\eta=(i,j)\in
M\setminus\Psi$ such that $i\in I$, $j\in J$, and $M_I^J$ is a
principal minor (i.e., the root $\eta$ lies inside a principal
minor).

First we show that there is a nonempty Zariski-open subset $V$ in
$\mathcal{Y}$ such that for any $x\in V$ there is a $h\in N$ such
that $\mathrm{Ad}_h\,x$ lies in (\ref{sum}) and any principal minor
in $\mathrm{Ad}_h\,x$ is different from zero.

Let $x\in\mathcal{Y}$. As above, we denote
$R_{k-1}=\displaystyle\sum_{t=1}^{k-1}r_t$ and
$\widetilde{r}=\widetilde{r}_k=\min(r,r_k)$. Let the principal minor
$M_I^J$ that stands to the right of the $k$th block in
$\mathfrak{r}$ be situated at the intersection of
$I=\{R_{k-1}+1,R_{k-1}+2,\ldots,R_{k-1}+\widetilde{r}\}$ rows and
$J=\{j_1,j_2,\ldots,j_{\widetilde{r}}\}$ columns.

By Lemma 3.2, there is a root from $\Psi$ in every column of the
principal minor, but not in any row of the principal minor there is
a root from $\Psi$. The roots $(4,7)$ and $(4,8)$ in Diagram 5 serve
as an example. These roots do not lie in the system $\Psi$, and the
fourth row belongs to the set of rows of the principal minors.

\medskip
\textbf{Example 3}. The diagram for a parabolic subalgebra with the
following sizes of diagonal block: $(1,1,4,2)$.
\begin{center}
{\begin{tabular}{|p{0.1cm}|p{0.1cm}|p{0.1cm}|p{0.1cm}|p{0.1cm}|
p{0.1cm}|p{0.1cm}|p{0.1cm}|c} \multicolumn{2}{l}{{\small
1\hspace{5pt} 2\!}}&\multicolumn{2}{l}{{\small 3\hspace{5pt}
4\!}}&\multicolumn{2}{l}{{\small 5\hspace{5pt} 6\!}}&
\multicolumn{2}{l}{{\small 7\hspace{5pt} 8\!}}\\
\cline{1-8} 1&$\boxtimes$&&&&&&&{\small1}\\
\cline{1-8} &1&$\boxtimes$&&&&&&{\small2}\\
\cline{2-8} \multicolumn{2}{|c|}{}&\multicolumn{4}{|l|}{1}&
$\boxtimes$&$\boxtimes$&{\small3}\\
\cline{7-8} \multicolumn{2}{|c|}{}&\multicolumn{4}{|c|}{1\quad\,}&
&&{\small4}\\
\cline{7-8}
\multicolumn{2}{|c|}{}&\multicolumn{4}{|c|}{\hspace{17pt}1}&
&$\otimes$&{\small5}\\
\cline{7-8}
\multicolumn{2}{|c|}{}&\multicolumn{4}{|r|}{1}&$\otimes$&
&{\small 6}\\
\cline{3-8} \multicolumn{6}{|c|}{}&\multicolumn{2}{|l|}{1}&{\small 7}\\
\multicolumn{6}{|c|}{}&\multicolumn{2}{|r|}{1}&{\small8}\\
\cline{1-8} \multicolumn{8}{c}{Diagram 5}\\
\end{tabular}}
\end{center}

Let $i\in I$, $j\in J$, and $(i,j)\not\in\Psi$. Then $(i,j)\not\in
S\cup\Phi$. Since $j$ is a column of the principal minor, there is a
root from $\Psi$ in the $j$th column. Therefore, a root
$(\widetilde{i},j)\in S$ exists; suppose $\widetilde{i}>i$. The
adjoint action on $x$ of the element
\begin{equation}
\widetilde{h}=E+t\cdot E_{i,\widetilde{i}},\ t\neq0,\label{XXX}
\end{equation}
enables one to obtain an arbitrary nonzero number at the entry
$(i,j)$ of the matrix $x$ ($w_{\widetilde{i},j}(x)\neq0$ by the
definition of $\mathcal{Y}$). We show that the adjoint
action~(\ref{XXX}) does not change the entries of the matrix $x$,
except for $(i,j)$. On the one hand, the number that stands at the
$(\widetilde{i},j)$ entry of the matrix~$x$ is a unique nonzero in
the $\widetilde{i}$th row. Indeed, if
$w_{\widetilde{i},\widetilde{j}}(x)\neq0$ for some~$\widetilde{j}$,
then $(\widetilde{i},\widetilde{j})\in\Phi$ and $\widetilde{j}<j$.
But in this case, a root $(b,\widetilde{i})\in S$ exists. Then, by
the corollary of Lemma 2.3, there is a root $(a,i)\in S$ in the
$i$th column.  Consequently, the pair of roots
$\left((a,i),(\widetilde{i},j)\right)$ is admissible, whence
$(i,j)\in\Phi$, which is impossible. On the other hand, to obtain
some changes in $x$ upon adjoint action~(\ref{XXX}), it is necessary
that the $i$th  column contain a root from the base $S$. Again we
have a contradiction with $(i,j)\not\in\Phi$. Thus we showed that
adjoint action (\ref{XXX}) does not change the matrix $x$, except
for the entry $(i,j)$.

Since the roots from $S$ occupy places not above the secondary
diagonal in a principal minor, there exists $h\in N$ such that
$\mathrm{Ad}_h\,x$ lies in (\ref{sum}) and
$w_{i,j}\left(\mathrm{Ad}_h\,x\right)\neq0$ for any root $(i,j)$
that stands not below the secondary diagonal in a principal minor.
then the set of $x\in\mathcal{Y}$ for which all principal minors in
$\mathrm{Ad}_h\,x$ are not equal to zero is a nonempty Zariski-open
set; denote it by $V$.

\medskip
Now we turn to a proof of the proposition. Denote
$x'=\mathrm{Ad}_h\,x$, and let $x'$ lie in (\ref{sum}). Suppose
$x\in V$. Then all principal minors in $x'$ are different from zero.

Consider the first $r_1$ rows in $x'$. Obviously, they satisfy
conditions 1--3 imposed on the matrices from $\mathcal{X}$, because
they do not contain roots from $\Phi$ and in any roe (column) there
is no more one root from $S$.

Now assume that for any row the number of which is no greater than
$R_{k-1}$, the matrix $x'$ satisfies requirements 1--3 in the
definition of the set $\mathcal{X}$. To prove the proposition, it is
sufficient to show that there exists a nonempty Zariski-open subset
$U\subset\mathcal{Y}$ such that for any $x\in U$ there is $g\in N$
such that the matrix $\mathrm{Ad}_g\,x'=\mathrm{Ad}_{gh}\,x$
satisfies requirements 1--3in the definition of the
set~$\mathcal{X}$ for rows the numbers of which do not exceed
$R_{k}$. We represent the part of the diagram that corresponds to
the rows $R_{k-1}+1,R_{k-1}+2,\ldots,R_k$. Lemma 3.2 implies that it
has the following form
\begin{center}
\begin{tabular}{|p{3.1cm}|p{0.2cm}|p{0.2cm}|p{0.2cm}|p{0.5cm}|
p{0.2cm}|p{0.2cm}|p{0.2cm}|p{0.2cm}|p{0.2cm}l}
\multicolumn{2}{r}{$j_1$}&\multicolumn{2}{r}{$j_2$}&
\multicolumn{2}{r}{...\ \ \ $j_{\widetilde{r}}$}\\
\cline{1-10} 1&$\boxtimes$&...&$\boxtimes$&$\ldots$&
$\boxtimes$&...&$\times$&&&$R_{k-1}+1$\\
\cline{2-10} \quad$\ldots$&\multicolumn{9}{|c}{$\dotfill$}\\
\cline{2-10} \qquad\ \ \!\!1&$\boxtimes$&...&$\boxtimes$&$\ldots$&
$\boxtimes$&...&$\otimes$&&&$R_{k-1}+\widetilde{r}-1$\\
\cline{2-10} \qquad\qquad
\!\!1&$\boxtimes$&...&$\boxtimes$&$\ldots$&
$\boxtimes$&&&&&$R_{k-1}+\widetilde{r}$\\
\cline{2-10} \qquad\qquad\ $\ldots$&\multicolumn{9}{|c}{$\dotfill$}\\
\cline{2-10} \qquad\qquad\qquad\!1&&...&$\otimes$&$\ldots$& & & & &\\
\cline{2-10} \qquad\qquad\qquad\quad 1&$\otimes$& & &$\ldots$& & & &
&&$R_{k}$\\
\cline{1-10}
\end{tabular}
\end{center}

Let the principal minor $M_I^J$ stand to the right of the block with
number $k$ in $\mathfrak{r}$ at the intersection of the rows
$I=\{R_{k-1}+1,R_{k-1}+2,\ldots,R_{k-1}+\widetilde{r}\}$ and the
columns $J=\{j_1,j_2,\ldots,j_{\widetilde{r}}\}$, where
$j_1<j_2<\ldots<j_{\widetilde{r}}$. Denote by $U'$ the nonempty
Zariski-open subset of those $x\in\mathcal{Y}$ for which the
principal minor $M_I^J$ in $x'$ is different from zero. We show that
for any $x\in U=U'\cap V$, there is $g\in N$ such that the matrix
$\mathrm{Ad}_g\,x'$ satisfies the requirements of the definition of
the set $\mathcal{X}$ for rows the numbers of which does not exceed
$R_k$.

Thus, let $x\in U$; therefore the principal minor $M_I^J$ is
different from zero. Then we can put it in triangular form with the
help of the adjoint action of element s from $N$ of the form
\begin{equation}
\widetilde{g}_{a}=E+t_1^a\cdot E_{j_1,j_a}+\ldots+t_{a-1}^a\cdot
E_{j_{a-1},j_a},\quad a=2,\ldots,\widetilde{r},\label{g1}
\end{equation}
choosing appropriate values of $t_i^a$, $i=1,\ldots,a-1$. In such a
way we achieve the fulfillment of condition 1 in the definition of
$\mathcal{X}$ for rows the number of which do not exceed $R_k$.

Next, let $(i,j)\in M\setminus\Psi$, where $i\in I$ and
$j>j_{\widetilde{r}}$, moreover $w_{i,j}(x')\neq0$. Similarly, we
act on $x'$ by conjugation by the element
\begin{equation}
\widehat{g}_{j}=E+t_1\cdot E_{j_1,j}+\ldots+t_{\widetilde{r}}\cdot
E_{j_{\widetilde{r}},j},
\end{equation}
where $t_1,\ldots,t_{\widetilde{r}}$ are appropriate values. We
achieve that the entry $(i,j)$ of the matrix
$\mathrm{Ad}_{\widehat{g}_{j}}\,x'$ becomes zero and the conjugation
by $\widehat{g}_j$ does not change the other entries of the matrix
in the $R_{k-1}+1,R_{k-1}+2,\ldots,R_{k}$ rows. The latter means
that condition 2 in the definition of the set $\mathcal{X}$ is
satisfied for the rows the numbers of which do not exceed $R_k$.

Now we assume that the principal minor to the right of the $k$th
block is triangular. Let $w_{i,j}(x')\neq0$ for a root $(i,j)$,
where $i\in I$, $j_{a-1}<j<j_a$ for some $j_{a-1},j_a\in J$, and the
root $(i,j_a)$ lies below the secondary diagonal in the principal
minor. It is sufficient to consider the conjugation by element
(\ref{g1}), replacing $j$ by $j_a$. Thus, the condition 3 in the
definition of $\mathcal{X}$ for rows the numbers of which do not
exceed $R_k$ is satisfied.~$\Box$

\medskip
Now let $\mathcal{S}$ be the set of denominators generated by minors
$M_\xi$, $\xi\in S$. We form the localization
$K[\mathfrak{m}]_\mathcal{S}$ of the algebra $K[\mathfrak{m}]$ on
$\mathcal{S}$.

Consider the Zariski-open subset in $\mathfrak{m}$
$$U_0=\left\{x\in\mathfrak{m}:M_{\xi}(x)\neq0\mbox{, для любого
}\xi\in S\right\}.$$ Then $K[\mathfrak{m}]_{\mathcal{S}}=K[U_0]$. We
define a regular map
\begin{equation}
\pi:U_0\rightarrow\mathcal{Y}\label{pi}
\end{equation}
by the rule
$$x\in U_0\mapsto \sum_{\xi\in S}\frac{M_{\xi}(x)}{\widehat{M}_{\xi}(x)}
\cdot E_{\xi}+\sum_{\varphi\in
\Phi}\frac{L_{\varphi}(x)}{M_{\xi_1}(x) \widehat{M}_{\xi_2}(x)}\cdot
E_{\varphi}\in\mathcal{Y},$$ where
\begin{itemize}
\item[1)] $\widehat{M}_{\gamma}(x)$ is a number $\prod
M_{\gamma'}(x)$, and the product proceeds over all roots $\gamma'$
such that $\gamma'<\gamma$ and $\gamma'$ is the greatest root in the
sense of the order~$<$ introduced above;
\item[2)] the root $\varphi\in\Phi$ from the second sum
corresponds to the admissible pair $q=(\xi_1,\xi_2)$.
\end{itemize}
Then $M_{\xi}(x)=M_{\xi}\big(\pi(x)\big)$,
$L_{\varphi}(x)=L_{\varphi}\big(\pi(x)\big)$ for any $\xi\in S$,
$\varphi\in\Phi$. Since $M_{\xi}$ and $L_{\varphi}$ are invariants,
we have $\pi\big(\mathrm{Ad}_g\,x\big)=\pi(x)$ for any $g\in N$.

\medskip
Now we are in a position to prove theorem 1.7.

\medskip
\textsc{Proof of Theorem 1.7} proceeds by induction on $n$.

For $n=2$ the statement is obvious. Assume that the theorem is valid
for matrices of sizes lesser than $n$. We prove it for matrices of
size $n$. Let the reductive part $\mathfrak{r}$ of a parabolic
subalgebra $\mathfrak{p}$ consists of blocks of sizes
$(r_1,r_2,\ldots,r_s)$, where $r_1+r_2+\ldots+r_s=n$. As above,
$r=r_s$. Introduce the notation
$$\mathfrak{z}_0=\left\{\left(\begin{tabular}{p{2.1cm}|p{0.3cm}}
&*\\ \quad\ \ \,{\Large0}&\ldots\\&*\\
\cline{1-2}
0\quad\ldots\quad0&0\\
\dotfill&\ldots\\
0\quad\ldots\quad0&0\\
\end{tabular}\right)\!\!\!\!\!\!\!\!
\raisebox{-22pt}{$\left.\begin{array}{c}\\
\\ \\\end{array}\right\}\!\emph{r}$}\,\right\},\quad
\mathcal{Z}_0=\mathrm{exp}\,\mathfrak{z}_0.$$

Since $\mathfrak{z}_0$ is an ideal in the subalgebra $\mathfrak{n}$,
we have $\mathrm{Ad}_N\,\mathfrak{z}_0=\mathfrak{z}_0$. We identify
$$\mathfrak{gl}(n-1,K)\mbox{ with the subalgebra }
\left\{\left(\begin{tabular}{p{2.1cm}|p{0.3cm}}
*\quad\ldots\quad*&0\\
\dotfill&\ldots\\
*\quad\ldots\quad*&0\\
\cline{1-2} 0\quad\ldots\quad0&0
\end{tabular}\right)\right\}\mbox{ in
}\mathfrak{gl}(n,K)$$ and
$$\mathrm{GL}(n-1,K)\mbox{ with the subgroup }
\left\{\left(\begin{tabular}{p{2.1cm}|p{0.3cm}}
*\quad\ldots\quad*&0\\
\dotfill&\ldots\\
*\quad\ldots\quad*&0\\
\cline{1-2} 0\quad\ldots\quad0&1
\end{tabular}\right)\right\}\mbox{ in }\mathrm{GL}(n,K).$$
Let $N_1$ be the subgroup of upper triangular matrices with unities
on the principal diagonal in $\mathrm{GL}(n-1,K)$,
$\mathfrak{n}_1=\mathrm{Lie\,N_1}$,
$\mathfrak{m}_1=\mathfrak{m}\cap\mathfrak{n}_1$, and let
$\mathcal{Y}_1$ be the natural projection of $\mathcal{Y}$ to
$\mathfrak{m}_1$.

Any element $x\in\mathfrak{m}$ is uniquely represented in the form
$x=x_1+z_1$, where $x_1\in\mathfrak{m}_1$ and
$z_1\in\mathfrak{z}_0$. By the inductive assumption, there exists a
nonempty Zariski-open subset $U_1\subset\mathfrak{m}_1$ such that
for any $x_1\in U_1$ there is $g'\in N_1$ such that
$\mathrm{Ad}_{g'}\,x_1\in\mathcal{Y}_1$. Since
$\mathfrak{m}=\mathfrak{m}_1\oplus\mathfrak{z}_0$, it follows that
$U_1+\mathfrak{z}_0$ is an open subset in $\mathfrak{m}$.

Proposition 4.3 implies the existence of a Zariski-open subset $V$
in $\mathcal{Y}$ such that for any $y\in V$, there is an element
$g_1\in N$ such that
$\widetilde{y}=\mathrm{Ad}_{g_1}\,y\in\mathcal{X}$. The preimage
$\pi^{-1}(V)$ of the set $V$ under map (\ref{pi}) is a Zariski-open
subset in $\mathfrak{m}$.

Denote
$$U=\left\{\begin{array}{ll}
U_1+\mathfrak{z}_0&\mbox{if there are no roots from }S\mbox{ in the
}n\mbox{th column};\\
\left(U_1+\mathfrak{z}_0\right)\cap\pi^{-1}(V)&\mbox{if there is a
root from
}S\mbox{ in the }n\mbox{th column}.\\
\end{array}\right.$$
From what we have said above, $U$ is a Zariski-open subset in
$\mathfrak{m}$.

\medskip
\textsc{Part I}. We show that for any $x\in U$ there exists $g\in N$
such that $\mathrm{Ad}_g\,x\in\mathcal{Y}$.

Thus, let $x\in U$. Denote $x'=\mathrm{Ad}_{g'}\,x$,
$y_0=\mathrm{Ad}_{g'}\,x_1$, $z_0=\mathrm{Ad}_{g'}\,z_1$. We obtain
$$x'=\mathrm{Ad}_{g'}\,x=\mathrm{Ad}_{g'}\,x_1+\mathrm{Ad}_{g'}\,
z_1=y_0+z_0,$$ где $y_0\in\mathcal{Y}_1$, $z_0\in\mathfrak{z}_0$.

\begin{itemize}
\item[1.] Consider the case where there are no roots from the base
$S$ in the $n$th column. In this case, $U=U_1+\mathfrak{z}_0$. We
show that there is $g\in N$ such that $\mathrm{Ad}_g\,x=y_0$. Lemma
2.1 implies that for any number $i$, where $1\leqslant i\leqslant
n-r$, there is a root $\xi\in S$ lying in the $i$th row (i.e.,
$\xi=(i,a_i)$ for some number $a_i$).  Since $y_0\in\mathcal{Y}_1$,
we have $w_{i,a_i}(x')=w_{i,a_i}(y_0)\neq0$. Denote
\begin{equation}
\begin{array}{rcl}
g_{n-r}=E+t_{n-r}\cdot E_{a_{n-r},n}\in\mathcal{Z}_0,&\mbox{where}&
t_{n-r}=\displaystyle\frac{w_{n-r,n}(z_0)}{w_{n-r,a_{n-r}}(y_0)};\\
g_i=E+t_i\cdot E_{a_i,n}\in\mathcal{Z}_0,&\mbox{where}&
t_i=\displaystyle\frac{w_{i,n}\big(\mathrm{Ad}_{g_{i+1}\ldots
g_{n-r}}\,z_0\big)}{w_{i,a_i}(y_0)}.\\
\end{array}
\end{equation}
Then
$$w_{a,n}\left(\mathrm{Ad}_{g_ig_{i+1}\ldots g_{n-r}}\, x'\right)=0$$
for any $a\geqslant i$. Denote
$$g=g_1\cdot g_2\cdot\ldots\cdot g_{n-r}\cdot g'.$$ We obtain
$$\begin{array}{c}
\mathrm{Ad}_g\,x=\mathrm{Ad}_{g_1g_2\ldots
g_{n-r}}\mathrm{Ad}_{g'}\,x= \mathrm{Ad}_{g_1g_2\ldots
g_{n-r}}\,x'=\\=\mathrm{Ad}_{g_1g_2\ldots
g_{n-r}}\,(y_0+z_0)=y_0\in\mathcal{Y},
\end{array}$$
which proves the statement in case 1.

\item[2.] Now consider the case where a root from
$S$ lies in the $n$th column. Here
$U=\left(U_1+\mathfrak{z}_0\right)\cap\pi^{-1}(V)$. Denote by
$(\widetilde{m},n)$ the root from $S$ lying in the $n$th column. Let
$(m,n)\in S\cup\Phi$ be a root which has the least number $m$.

Lemma 2.1 implies that for every $(i,n)\in M$, where
$i>\widetilde{m}$, a root from $S$ exists in the $i$th row.
Similarly to the first item, there is an element $h\in\mathcal{Z}_0$
such that $\mathrm{Ad}_h\,x'$ has zeros in the last column below the
root $(\widetilde{m},n)$ (i.e.,
$w_{a,n}\left(\mathrm{Ad}_h\,x'\right)=0$ for any
$a>\widetilde{m}$).

Now let $(i,n)\not\in S\cup\Phi$ and $m<i<\widetilde{m}$. Similarly
to the first part of the proof of Proposition 4.3, acting on $x'$ by
conjugation by element~(\ref{XXX}) for
\raisebox{0pt}[13pt]{$\widetilde{i}=\widetilde{m}$}, we can obtain
$$w_{i,n}(\mathrm{Ad}_h\,x')=0\mbox{ and
}w_{a,b}(\mathrm{Ad}_h\,x')=w_{a,b}(x')$$ for any root
$(a,b)\neq(i,n)$.

In view of the above-said, in the sequel we assume that for $x'$,
$w_{i,n}(x')=0$ for any root $(i,n)\not\in S\cup\Phi$ such that
$i>m$.

Consider the subspace $\mathfrak{z}\subset\mathfrak{z}_0$ that
consists of $z\in\mathfrak{z}_0$ for which $w_{i,n}(z)=0$ for
$i\geqslant m$. Denote by $z$ the natural projection of $z_0$ to
$\mathfrak{z}$. Then the matrix $y=x'-z$ lies in $\mathcal{Y}$.
Therefore, $x'=y+z$, $y\in\mathcal{Y}$, $z\in\mathfrak{z}$. From the
definition of invariants, $y=\pi(x')$. Since $x'\in\pi^{-1}(V)$, we
have $y\in V$.

From Proposition 4.3 it follows that there exists an element $g_1\in
N$ such that $\widetilde{y}=\mathrm{Ad}_{g_1}\,y\in\mathcal{X}$.
Then $\mathrm{Ad}_{g_1}\,x'=\widetilde{y}+\mathrm{Ad}_{g_1}\,z$.
Since $\mathfrak{z}$ is an ideal of $\mathfrak{n}$, we have
$\mathrm{Ad}_{g_1}\,z\in\mathfrak{z}$.

We show that there exists $g_2\in\mathcal{Z}_0$ such that
$$\mathrm{Ad}_{g_2g_1}\,x'=\widetilde{y}+\widetilde{z},$$
where $\widetilde{y}\in\mathcal{X}$ and
$\widetilde{z}\in\mathcal{T}\cap\mathfrak{z}$. Let
$\zeta=(i,n)\not\in T$, $i<m$. In the $n$th column, there is a root
from $S$ and thus from $\Psi$. Then, by definition of~$T$, there
exists a principal minor $M_I^J$, where $I$ is a set of numbers of
rows and $J$ is the set of numbers of columns, such that $i\in I$.
Let $\widetilde{r}$ be its size and
$J=\{j_1,\ldots,j_{\widetilde{r}}\}$,
$j_1<\ldots<j_{\widetilde{r}}$. Since $\widetilde{y}\in\mathcal{X}$,
the principal minor $M_I^J$ is different from zero. We choose
numbers $t_1,\ldots,t_{\widetilde{r}}$ in such a way that for the
element
$$g_{\zeta}=E+t_1\cdot E_{j_1,n}+\ldots+t_{\widetilde{r}}\cdot
E_{j_{\widetilde{r}},n}$$ the following conditions are satisfied:
\begin{itemize}
\item[(a)] $\mathrm{Ad}_{g_{\zeta}}\left(\widetilde{y}+
\mathrm{Ad}_{g_1}\,z\right)=\widetilde{y}+z'$, where
$z'\in\mathfrak{z}$;
\item[(b)] $w_{i,n}\left(\mathrm{Ad}_{g_{\zeta}g_1}\,z\right)=0$;
\item[(c)] if $i\neq j$, then the adjoint action of
$g_{\zeta}$ does not change the entry $(j,n)$ of the matrix
$\widetilde{y}+\mathrm{Ad}_{g_1}\,z$.
\end{itemize}
Then for $g_2=\displaystyle\prod_{\zeta=(i,n)\not\in
T}g_{\zeta}\in\mathcal{Z}_0$, we have
$$\mathrm{Ad}_{g_2g_1}\,x'=\widetilde{y}+\widetilde{z}\mbox{, where }
\widetilde{y}\in\mathcal{X}\mbox{ and }
\widetilde{z}\in\mathcal{T}\cap\mathfrak{z}.$$

Next, Proposition 4.2 implies the existence of $g_3\in N$ such that
$$\mathrm{Ad}_{g_3}\left(\widetilde{y}+\widetilde{z}\right)=
\widetilde{y}.$$ We obtain
$$\mathrm{Ad}_{g_3g_2g_1}\,x'=\mathrm{Ad}_{g_3g_2g_1}(y+z)=
\mathrm{Ad}_{g_3}(\widetilde{y}+\widetilde{z})=\widetilde{y}.$$ It
remains to apply the inverse action $g_1^{-1}$:
$$\mathrm{Ad}_{g_1^{-1}g_3g_2g_1g'}\,x=
\mathrm{Ad}_{g_1^{-1}g_3g_2g_1}\,x'=
\mathrm{Ad}_{g_1^{-1}}\,\widetilde{y}=
\mathrm{Ad}_{g_1^{-1}}\left(\mathrm{Ad}_{g_1}\,y\right)=
y\in\mathcal{Y}.$$
\end{itemize}

\medskip
\textsc{Part II}. Uniqueness follows from the fact that the map
(\ref{pi}) is constant on an $N$-orbit.~$\Box$

\medskip
Consider the localization $K[\mathfrak{m}]^N_\mathcal{S}$ of the
algebra of invariants $K[\mathfrak{m}]^N$ with respect to
$\mathcal{S}$. Since the minors $M_\xi$ are $N$-invariants, we have
$$K[\mathfrak{m}]^N_\mathcal{S}=\big(K[\mathfrak{m}]_S\big)^N.$$

\medskip
The following theorem was proved in~\cite{PS}; we repeat its proof
to maintain the understanding of the matter.

\medskip
\textbf{Theorem 4.4.} \emph{The ring $K[\mathfrak{m}]^N_S$ is the
ring of polynomials in} $M_\xi^{\pm 1}$, $\xi\in S$, \emph{and in}
$L_{\varphi}$, $\varphi\in\Phi$.

\medskip
\textsc{Proof.} Consider the restriction homomorphism
$\varpi:f\mapsto f|_\mathcal{Y}$ of the algebra $K[\mathfrak{m}]^N$
to $K[\mathcal{Y}]$. The image $\varpi(M_\xi)$ is equal to the
product
$$\pm x_\xi x_{\xi_1}\ldots x_{\xi_t},$$
where every $\xi_{i}$ is less than $\xi$ in the sense of above
order~$>$. We extend $\varpi$ to a homomorphism
$$\varpi_S:\,K[\mathfrak{m}]^N_\mathcal{S}\to K[\mathcal{Y}]_S,$$
where $K[\mathcal{Y}]_S$ is a localization of $K[\mathcal{Y}]$ with
respect to $x_\xi$, $\xi\in S$. We show that $\varpi_S$ is a
isomorphism.

If $f\in\mathrm{Ker}\,p_S$, then $f(\mathrm{Ad}_N\,\mathcal{Y})=0$.
Since, by Theorem 1.7, $\mathrm{Ad}_N\,\mathcal{Y}$ contains a
Zariski-open subset, we have $f=0$. Consequently, $\varpi_S$ is an
embedding of $K[\mathfrak{m}]^N_\mathcal{S}$ in $K[\mathcal{Y}]_S$.
Next, we write the formulas in the form
$$\begin{array}{l}
\varpi\left(M_\xi\right)=\pm x_\xi
\varpi\left(\widehat{M}_{\xi}\right),\\
\varpi\left(L_\varphi\right)=\pm x_\varphi
\varpi\left(M_{\gamma}\right)\varpi\left(
\widehat{M}_{\gamma'}\right),\end{array}$$ where the root
$\varphi\in\Phi$ corresponds to the admissible pair
$(\gamma,\gamma')$, and $\widehat{M}_{\xi}$ (
$\widehat{M}_{\gamma'}$, respectively) is the product $\prod
M_{\xi_1}$ ($\prod M_{\gamma_1}$, respectively), and the product
proceeds over the greatest roots $\xi_1$ such that $\xi_1$ is less
than $\xi$ ($\gamma'$, respectively) in the sense of the order~$<$.
We get
\begin{equation}
\begin{array}{l}\varpi_S\left(M_\xi
\widehat{M}^{-1}_{\xi}\right)=\pm x_\xi,\\
\varpi_S\left(L_\varphi
M^{-1}_{\gamma}\widehat{M}^{-1}_{\gamma'}\right)=\pm
x_\varphi.\end{array}\label{pi_S}
\end{equation}
From (\ref{pi_S}) it follows that the image of $\varpi_S$ coincides
with $K[\mathcal{Y}]_S$. Thus, $\varpi_S$ is an isomorphism.~$\Box$

\medskip
We obtain Theorem 1.8 as a consequence of Theorem 4.4.

\textsc{Department of Mechanics and Mathematics, Samara State
University, Russia}\\ \emph{E-mail address}:
\verb"victoria.sevostyanova@gmail.com"


\begin{thebibliography}{10}

\bibitem[B]{B}
M. Brion, Representations exceptionnelles des groups semi-simple,
\emph{Ann. Scient. Ec. Norm. Sup.} \textbf{18} (1985), 345--387.
\bibitem[GG]{GG}
М. Goto and F. Grosshans, Semisimple Lie algebras, Lect. Notes in
Pure Appl. Math., vol. 38, 1978.
\bibitem[K]{K}
H. Kraft, Geometrische Methoden in der Invariantentheorie, Friedr.
Vieweg and Sohn, Braunschweig/Wiesbaden (1985).
\bibitem[PS]{PS}
A. N. Panov and V. V. Sevostyanova, Regular $N$-orbits in the
nilradical of a parabolic subalgebra, \emph{Vestnik SamGU},
\textbf{7}(57) (2007), 152--161. See also
http://arxiv.org/abs/1203.2754.
\bibitem[PV]{PV}
V. L. Popov and E. B. Vinberg, Invariant theory, in:
\emph{Progressin Science and Technology}, VINITI, Moscow (1989), pp.
137--309.
\bibitem[R]{R}
R. W. Richardson, Conjugacy classes in parabolic subgroups of
semisimple  algebraic groups, \emph{Bull. London Math. Soc.}
\textbf{6} (1974), \mbox{21--24.}

\end{thebibliography}
\end{document}